\documentclass{birkjour}
\usepackage{hyperref}

\newtheorem{theorem}{Theorem}[section]
\theoremstyle{plain}

\newtheorem{definition}{Definition}

\newtheorem{lemma}{Lemma}[section]

\numberwithin{equation}{section}

\begin{document}

%
%
%
%
%
%
%
%
%

\title[New integral inequalities via $P$-convexity]{New integral inequalities via $P$-convexity}

\author[W. J. Liu]{Wenjun Liu}

\address{College of Mathematics and Statistics\\
Nanjing University of Information Science and Technology \\
Nanjing 210044, China}

\email{wjliu@nuist.edu.cn}

\thanks{This work was partly supported by the National Natural Science Foundation
of China (Grant No. 40975002) and the Natural Science Foundation of the Jiangsu
Higher Education Institutions (Grant No. 09KJB110005).}
\subjclass{Primary 26D15; Secondary 33B15, 26D07}

\keywords{Hermite's  inequality, H\"{o}lder's inequality, integral inequality,  $P$-convexity}

\date{January 28, 2012}

\begin{abstract}In this note we extend some new estimates of the integral $\int_a^b (x-a)^p(b-x)^qf(x)dx$ for functions when a power of the absolute value is $P-$convex.
\end{abstract}

\maketitle
\section{INTRODUCTION}

Let $I$ be an interval in $\mathbb{R}$. Then $f:I\rightarrow \mathbb{R}$ is said to be convex
if
\begin{equation*}
f\left( t x+\left( 1-t \right) y\right) \leq t f\left(
x\right) +\left( 1-t \right) f\left( y\right)
\end{equation*}%
holds for all $x,y\in I$ and $t \in \left[ 0,1\right]$.

The notion of quasi-convex functions generalizes the notion
of convex functions. More precisely, a function $f:[a,b]\rightarrow \mathbb{R}$ is said to be quasi-convex on $[a,b]$ if%
\begin{equation*}
f(t x+(1-t )y)\leq \max \{f(x),f(y)\}
\end{equation*}%
holds for any $x,y\in \lbrack a,b]$ and $t \in \lbrack 0,1]$. Clearly, any
convex function is a quasi-convex function. Furthermore, there exist
quasi-convex functions which are not convex (see \cite{daI}).

The generalized quadrature formula of Gauss-Jacobi type has the form
\begin{equation}\label{1.2'}
\int_a^b (x-a)^p(b-x)^qf(x)dx=\sum\limits_{k=0}^m B_{m,k}f(\gamma_k)+\mathcal{R}_m[f]
\end{equation}
for certain $B_{m,k}, \gamma_k$ and rest term $\mathcal{R}_m[f]$ (see \cite{SCB}).

In \cite{OSA},  \"{O}zdemir et al.  established several integral inequalities concerning the left-hand side of \eqref{1.2'} via some kinds of convexity.
Especially, they discussed the following
result connecting with  quasi-convex function:
\begin{theorem}\label{th1.1}
Let $f: [a,b]\rightarrow\mathbb{R}$ be continuous on $[a,b]$ such that $f \in L([a,b])$, $0\le a<b<\infty$. If $f$ is quasi-convex on $[a,b]$, then for some fixed   $p, q > 0$, we have
\begin{align*}
 \int_a^b (x-a)^p(b-x)^qf(x)dx
\leq (b-a)^{p+q+1} \beta(p+1, q+1)\max\{f(a), f(b)\},
\end{align*}
where $\beta (x, y)$ is the Euler Beta function.
\end{theorem}

Recently, Liu  \cite{l2012} established some new integral inequalities for quasi-convex functions as follows:

\begin{theorem}\label{th1.2}
Let $f: [a,b]\rightarrow\mathbb{R}$ be continuous on $[a,b]$ such that $f \in L([a,b])$, $0\le a<b<\infty$ and let $k>1$. If $|f|^{\frac{k}{k-1}}$ is quasi-convex on $[a,b]$,  for some fixed  $p, q > 0$, then
\begin{align*}
 &\int_a^b (x-a)^p(b-x)^qf(x)dx\\
\leq   & (b-a)^{p+q+1} \left[\beta(kp+1, kq+1)\right]^{\frac{1}{k}}\left( \max \left\{ \left| f(a)\right| ^{\frac{k}{k-1}},
\left| f(b)\right| ^{\frac{k}{k-1}}\right\} \right) ^{\frac{k-1}{k}}.
\end{align*}
\end{theorem}

\begin{theorem}\label{th1.3}
Let $f: [a,b]\rightarrow\mathbb{R}$ be continuous on $[a,b]$ such that $f \in L([a,b])$, $0\le a<b<\infty$ and let $l\ge 1$. If $|f|^{l}$ is quasi-convex on $[a,b]$,  for some fixed  $p, q > 0$, then
\begin{align*}
 &\int_a^b (x-a)^p(b-x)^qf(x)dx\\
\leq  &(b-a)^{p+q+1} \beta(p+1, q+1) \left( \max \left\{ \left| f(a)\right| ^l,
\left| f(b)\right| ^l\right\} \right) ^{\frac{1}{l}}.
\end{align*}
\end{theorem}

On the other hand, Dragomir et al. in \cite{dpp1995} defined the following class of functions.

\begin{definition}
\label{d1} Let $I\subseteq R$ be an interval. The function $f: \rightarrow\mathbb{R}$ is said to belong to
the class $P(I)$ (or to be $P$-convex) if it is nonnegative and, for all $x, y \in I$ and $t\in [0, 1]$,
satisfies the inequality
\begin{equation*}
f(t x+(1-t )y)\leq  f(x)+f(y).
\end{equation*}%
\end{definition}

Note that $P(I)$ contain all nonnegative convex and quasiconvex functions. Since
then numerous articles have appeared in the literature reflecting further applications in
this category; see \cite{ADD, ao2010, Dragomir, DP, dp1998, HC, oy2011, SSO, Set2, SSOR, tyd2003} and references therein.

The main purpose of this note is to establish some new estimates of the integral $\int_a^b (x-a)^p(b-x)^qf(x)dx$ for functions when a power of the absolute value is $P-$convex. That is, this study
is a further continuation and generalization of \cite{OSA} via $P$-convexity.

\section{New integral inequalities via $P$-convexity}\label{se3}

In this section we generalize Theorems \ref{th1.1}-\ref{th1.3} with a $P$-convex function setting. For this purpose, we need the
following lemma (see \cite[Lemma 2.1]{OSA}):

\begin{lemma}\label{le1.1} Let $f: [a,b]\subset [0, \infty)\rightarrow\mathbb{R}$ be continuous on $[a,b]$ such that $f \in L([a,b])$, $a<b$. Then the equality
\begin{equation}\label{1.2}
\int_a^b (x-a)^p(b-x)^qf(x)dx=(b-a)^{p+q+1} \int_0^1 (1-t)^pt^qf(ta+(1-t)b)dt
\end{equation}
holds for some fixed $p, q > 0$.
\end{lemma}

The next theorem gives a new result for $P$-convex functions.

\begin{theorem}\label{th3.0}
Let $f: [a,b]\rightarrow\mathbb{R}$ be continuous on $[a,b]$ such that $f \in L([a,b])$, $0\le a<b<\infty$. If $|f|$ is $P$-convex on $[a,b]$,  for some fixed  $p, q > 0$, then
\begin{align}\label{3.0}
  &\int_a^b (x-a)^p(b-x)^qf(x)dx\nonumber\\
\leq  & (b-a)^{p+q+1} \beta(p+1, q+1) (|f(a)| +  |f(b)|),
\end{align}
where $\beta (x, y)$ is the Euler Beta function.
\end{theorem}

\begin{proof}
By Lemma \ref{le1.1}, the Beta function which is defined for $x, y > 0$ as
$$\beta (x, y)=\int_0^1 t^{x-1}(1-t)^{y-1}dt$$
 and the fact that $f$ is $P$-convex  on $[a,b]$, we have
\begin{align*}
&\int_a^b (x-a)^p(b-x)^qf(x)dx\\
\leq &(b-a)^{p+q+1} \int_0^1   (1-t)^pt^q  |f(ta+(1-t)b)|dt\\
\leq & (b-a)^{p+q+1}
 \int_0^1 (1-t)^{p}t^{q}  (|f(a)|+|f(b)|) dt  \\
= & (b-a)^{p+q+1}   \beta(q+1, p+1)(|f(a)| +  |f(b)| ),
\end{align*}
which completes the proof.
\end{proof}

The corresponding version for powers of the absolute value is
incorporated in the following result.

\begin{theorem}\label{th3.1}
Let $f: [a,b]\rightarrow\mathbb{R}$ be continuous on $[a,b]$ such that $f \in L([a,b])$, $0\le a<b<\infty$ and let $k>1$. If $|f|^{\frac{k}{k-1}}$ is $P$-convex on $[a,b]$,  for some fixed  $p, q > 0$, then
\begin{align}\label{3.1}
 &\int_a^b (x-a)^p(b-x)^qf(x)dx\nonumber\\
\leq   & (b-a)^{p+q+1} \left[\beta(kp+1, kq+1)\right]^{\frac{1}{k}}\left(   \left| f(a)\right| ^{\frac{k}{k-1}}+
\left| f(b)\right| ^{\frac{k}{k-1}}  \right) ^{\frac{k-1}{k}}.
\end{align}
\end{theorem}

\begin{proof}
By Lemma \ref{le1.1}, H\"older's inequality, the definition of Beta function and the fact that $|f|^{\frac{k}{k-1}}$ is $P$-convex on $[a,b]$, we have
\begin{align*}
&\int_a^b (x-a)^p(b-x)^qf(x)dx\\
\leq & (b-a)^{p+q+1}\left[\int_0^1 (1-t)^{kp}t^{kq}dt\right]^\frac{1}{k}
\left[\int_0^1 |f(ta+(1-t)b)|^{\frac{k}{k-1}} dt\right]^\frac{k-1}{k} \\
\leq & (b-a)^{p+q+1}\left[\beta(kq+1, kp+1)\right]^{\frac{1}{k}} \left[\int_0^1   \left( \left| f(a)\right| ^{\frac{k}{k-1}}+
\left| f(b)\right| ^{\frac{k}{k-1}}\right)dt \right]^\frac{k-1}{k}\\
= & (b-a)^{p+q+1}\left[\beta(kq+1, kp+1)\right]^{\frac{1}{k}} \left(    \left| f(a)\right| ^{\frac{k}{k-1}}+
\left| f(b)\right| ^{\frac{k}{k-1}}   \right)^\frac{k-1}{k},
\end{align*}
which completes the proof.
\end{proof}

A more general inequality using Lemma \ref{le1.1} is as follows:

\begin{theorem}\label{th3.2}
Let $f: [a,b]\rightarrow\mathbb{R}$ be continuous on $[a,b]$ such that $f \in L([a,b])$, $0\le a<b<\infty$ and let $l> 1$. If $|f|^{l}$ is $P$-convex on $[a,b]$,  for some fixed  $p, q > 0$, then
\begin{align}\label{3.2}
 &\int_a^b (x-a)^p(b-x)^qf(x)dx\nonumber\\
\leq & (b-a)^{p+q+1} \beta(p+1, q+1) \left(   \left| f(a)\right| ^l+
\left| f(b)\right| ^l \right) ^{\frac{1}{l}},
\end{align}
where $\beta (x, y)$ is the Euler Beta function.
\end{theorem}

\begin{proof}
By Lemma \ref{le1.1}, H\"older's inequality, the definition of Beta function and the fact that $|f|^{l}$ is $P$-convex  on $[a,b]$, we have
\begin{align*}
&\int_a^b (x-a)^p(b-x)^qf(x)dx\\
=&(b-a)^{p+q+1} \int_0^1 \left[(1-t)^pt^q\right]^\frac{l-1}{l} \left[(1-t)^pt^q\right]^\frac{1}{l} f(ta+(1-t)b)dt\\
\leq & (b-a)^{p+q+1}\left[\int_0^1 (1-t)^{p}t^{q}dt\right]^\frac{l-1}{l}
\left[\int_0^1 (1-t)^{p}t^{q} |f(ta+(1-t)b)|^{l} dt\right]^\frac{1}{l} \\
\leq & (b-a)^{p+q+1}\left[\beta(q+1, p+1)\right]^{\frac{l-1}{l}}  \left[ \left( \left| f(a)\right| ^l+
\left| f(b)\right| ^l\right) \beta(q+1, p+1) \right]^\frac{1}{l}\\
=&(b-a)^{p+q+1} \beta(p+1, q+1) \left(   \left| f(a)\right| ^l+
\left| f(b)\right| ^l \right) ^{\frac{1}{l}},
\end{align*}
which completes the proof.
\end{proof}

\end{document}